\documentclass[letterpaper, 10 pt, conference]{ieeeconf}  

\IEEEoverridecommandlockouts                              

\usepackage{times} 
\usepackage{amsmath} 
\usepackage{amssymb}  
\usepackage{siunitx}

\usepackage{xcolor}
\usepackage{pgfplots}             

\usepackage{tikz}             
\usepackage{verbatim}
\usepackage{nicefrac}
\usepackage{mathtools}

\usetikzlibrary{pgfplots.groupplots}

\title{\LARGE \bf
Algebraic Fault Detection and Identification for Rigid Robots
}

\author{Alexander Lomakin and Joachim Deutscher
\thanks{Alexander Lomakin and Joachim Deutscher are with the Lehrstuhl f{\"u}r Regelungstechnik,
        Friedrich-Alexander-Universit{\"a}t, Erlangen-N{\"u}rnberg, Germany
        {\tt\small $\{$alexander.lomakin, joachim.deutscher$\}$@fau.de}}
}

\begin{document}

\maketitle
\thispagestyle{empty}
\pagestyle{empty}

\begin{abstract}
 This paper presents a method for algebraic fault detection and identification of nonlinear mechanical systems, describing rigid robots, by using an approximation with orthonormal Jacobi polynomials. An explicit expression is derived for the fault from the equation of motion, which is decoupled from disturbances and only depends on measurable signals and their time derivatives. Fault detection and identification is then achieved by polynomial approximation of the determined fault term. The results are illustrated for a faulty SCARA.
\end{abstract}

\section{INTRODUCTION}
	Due to increasing demands on mechanical systems, e.g. robots, in terms of availability and safety and the simultaneously increasing complexity, there is a continuous need of detecting and identifying potentially occurring faults in real time such that the general safety of the entire system can be guaranteed. This is why fault detection and identification methods are also becoming more and more important in industrial applications. Based on model knowledge of the system and an FMEA (Failure Mode and Effects Analysis), faults as well as other disturbances can be considered as unknown inputs in the motion equations. This fact is often used for fault detection and identification, especially for the observer-based approaches (see e.g. \cite{Chen:1999:RMF:316586}), whereas these methods are well established for linear systems (see \cite{Ding2008.}). Problems with these methods, however, lie in the selection of residuals and in the high implementation and parameterization effort for the mostly nonlinear systems. Furthermore, the high effort required for disturbance decoupling (see e.g. \cite{Persis2001AGA}) and the nonlinear observer design render these approaches uninteresting for practical applications. With algebraic approaches (see e.g. \cite{mai:inria-00159308}) these problems are avoided. By reconstructing the time derivatives of the input and output variables with the help of derivative estimators the fault can be determined algebraically from the equation of motion.        

 In this paper an algebraic method for fault detection and identification is introduced. It is based on an algebraic polynomial approximation of the fault by Jacobi polynomials (see \cite{6502453}). The properties of the approximation are introduced and then used to successively eliminate the time derivatives such that the fault reconstruction can be determined solely from the existing measurement signals, independent of the characteristics of the input and output signals and the disturbance.

In the next section a formulation of the considered fault detection problem is given. Then, the polynomial approximation and its properties are presented in Section 3. In Sections 4 and 5 these results are used for a general fault detection and fault identification. The proposed method is demonstrated by means of a faulty SCARA.

\section{PROBLEM FORMULATION}

Consider a robot as a general nonlinear mechanical system, which has $n$ fully actuated rigid joints. The motion of the robot can be described by the generalized coordinates $q\in \mathbb{R}^{n}$, such as link positions that are measured and whose associated time derivatives are $\dot{q}\in \mathbb{R}^{n}$ and $\ddot{q} \in \mathbb{R}^{n}$. Taking the fault $f$ and the disturbance $d$ into account, the dynamic behaviour of the robot can be described by
\begin{subequations} \label{generalMotion}
	\begin{align}
		\mathcal{M}(q)\,\ddot{q} +\mathcal{C}\left(q,\dot{q}\right)\dot{q}+\mathcal{G}\left(q \right) &= \,u + \mathcal{F}(q,\dot{q})f \nonumber\\
 		& \quad +  \mathcal{D}(q,\dot{q})d \label{generalMotionA}\\
	y &= q \label{generalMotionOutput}
	\end{align}
\end{subequations} 
with the initial condition $q(0), \dot{q}(0) \in \mathbb{R}^{n}$. In this formulation, $\mathcal{M}(q) \in \mathbb{R}^{n\times n}$ corresponds to the generalized inertia matrix and contains mass, inertia and geometric parameters of the joints. The vector $\mathcal{C}\left(q,\dot{q}\right)\dot{q} \in \mathbb{R}^{n}$ represents the Coriolis and centripetal components and $\mathcal{G}\left(q \right)\in \mathbb{R}^{n}$ is regarded as the influence of gravitational force. The output $y \in \mathbb{R}^{n}$ of \eqref{generalMotion} is available for measurement and equal to the link positions $q$. The presented method implicitly reduces the influence of measurement noise (see, e.g. \cite{Mboup2009}). Hence, the influence of noise on the system is not discussed in this paper. The right-hand side of \eqref{generalMotion} contains all the non-conservative forces that affect the robot. These forces include the input torque $u\in \mathbb{R}^{n}$, as well as the unknown fault $f= \mathrm{col}(f_1,f_2,...,f_{n_f}) \in \mathbb{R}^{n_f}$ and unknown disturbance $d =\mathrm{col}(d_1,d_2,...,d_{n_d}) \in \mathbb{R}^{n_d}$. The corresponding matrices 
\begin{equation}
\mathcal{F}(q,\dot{q}) = \bigl[\mathcal{F}_1(q,\dot{q}) \, ...  \, \mathcal{F}_{n_f}(q,\dot{q})\bigr] \in \mathbb{R}^{n \times n_f}
\end{equation}
 and 
 \begin{equation}
 \mathcal{D}(q,\dot{q}) = \bigl[\mathcal{D}_1(q,\dot{q}) \,  ...  \, \mathcal{D}_{n_d}(q,\dot{q})\bigr] \in \mathbb{R}^{n \times n_d}
 \end{equation} 
 represent the influence of the fault and the disturbance, respectively, and are assumed to be known. It is also assumed that $\mathcal{D}(q,\dot{q})$ has full column rank and therefore $\mathrm{rank}\,\mathcal{D}(q,\dot{q}) = n_d$. In this paper, faults are assumed to be additive only, e.g., actuator faults. 

The fault detection and identification problem addressed in this paper is interpreted analogously to \cite{isermann2005fault} as the subsequent design problem.

For a given system \eqref{generalMotion}, find a \textit{residual signal} $r$ 
\begin{equation}
\label{residual1}
	r = \Phi(q,\dot{q},\ddot{q},u)  \in \mathbb{R}^{n_r}
\end{equation}
such that the conditions
\begin{enumerate}
	\item[I.] $\lim\limits_{t\rightarrow \infty}\lVert r\rVert = 0, \qquad  \forall \, f= 0$  
	\item[II.] $ r \ne 0, \qquad \qquad \quad  \forall \; f\ne 0$  
\end{enumerate}
are satisfied for any input $u \in \mathbb{R}^n$, any disturbance $d\in \mathbb{R}^{n_d}$ and any initial states $q(0), \dot{q}(0) \in \mathbb{R}^n$. Then the fault $f$ can be detected by the residual signal $r$. If additionally to I and II  any two particular faults $f_i$ and $f_j$, with $ f_i,f_j \in  \mathbb{R}^{n_f}$ and $f_j \ne f_i $, the corresponding residual signals $r_i$ and $r_j$ can be distinguished for any finite time interval $\mathcal{I}_{t}$, the fault $f$ can also be isolated by the residual signal $r$. Furthermore, if $f_i$ and $f_j$ can be additionally estimated by $r$ for any finite time interval $\mathcal{I}_{t}$, the fault is assumed to be identifiable. Henceforth, the addressed problem will be regarded as the residual generation problem for fault detection (RGP-FD) and fault identification (RGP-FDI), respectively.

\section{ALGEBRAIC FAULT DETECTION AND IDENTIFICATION}
To solve the presented fault detection problem, it is necessary to calculate the residual signal $r$. Since the derivatives of $q$ are not available for measurement, a computational method to determine the time derivatives $\dot{q},\ddot{q}$ and thus the residual, which only needs the available signals $u$ and $q$, will be presented in this section. 
\subsection{Polynomial approximation} 
As already shown in \cite{6502453} and \cite{DBLP:conf/cdc/KiltzR13}, a polynomial approximation of a function can be determined to calculate the derivatives without differentiation. In this paper an analogous method is used by defining a polynomial approximation operator $\mathcal{P}\{\cdot\}$, which calculates the polynomial approximation of a function by an integral transformation. Therefore, the derivation and definition of $\mathcal{P}\{\cdot\}$ will be described first.

Consider the function $x\in L_2(\left[t-T,t\right])$ defined on the sliding time window $\mathcal{I}_{t,T} = \left[ t-T, t\right]$, $T>0$. In order to only take a constant time independent approximation interval $\tilde{\mathcal{I}} = [-1, 1]$ into account, the bijective transformation $\phi_T: \tilde{\mathcal{I}}  \mapsto  \mathcal{I}_{t,T}$, which maps the interval $\tilde{\mathcal{I}} = [-1, 1]$ to the given time window $\mathcal{I}_{t,T}$ is introduced. This transformation reads
 \begin{equation}
 \label{timeTransform1}
\phi_T(\tau') = t + \frac{\tau'-1}{2}T, \qquad \tau' \in \left[-1,1\right],
 \end{equation}
 with the inverse mapping $\phi^{-1}_T:   \mathcal{I}_{t,T} \mapsto  \tilde{\mathcal{I}}$ given by 
 \begin{equation}
 \label{timeTransform2}
 \phi^{-1}_T(t') = 1 + 2\frac{t'-t}{T}, \qquad t' \in \left[t-T,t\right].
 \end{equation}
 The transformed function $\bar{x} = x \circ \phi_T $ is then defined on a Hilbert space $\mathcal{H} = L_2(\left[-1,1\right])$ with the inner product 
 \begin{equation}
 \label{scalar}
 \langle \varphi_i,\varphi_j \rangle = \int_{-1}^{1} \varphi_i(\tau) \varphi_j(\tau) w^{(\alpha, \beta)}(\tau) \mathrm{d}\tau, \quad \forall \varphi_i,\varphi_j \in \mathcal{H},
 \end{equation}
  and the induced norm 
 \begin{equation}
 \label{norm}
 \left\lVert \varphi \right\rVert = \sqrt{\langle \varphi,\varphi \rangle}, \quad \forall \varphi \in \mathcal{H}.
 \end{equation}
The weight function $w^{(\alpha, \beta)}$, which allows to consider Jacobi polynomials as an orthonormal basis for $\mathcal{H}$, is given by 
 \begin{equation}
 \label{weight}
 w^{(\alpha, \beta)}(\tau) = \begin{cases}
 (1-\tau )^{\alpha } (1+\tau)^{\beta}, &\tau \in \left[-1,1\right],\\
 0, & \tau \notin \left[-1,1\right],
 \end{cases} 
 \end{equation}
with the real exponential coefficients $\alpha, \beta > -1$ as a degree of freedom. They can be chosen to achieve special blocking properties of individual frequencies in the signal by approximation (see, e.g., \cite{DBLP:conf/cdc/KiltzR13}). 
Then, it is possible to introduce an orthonormal basis $\{P_i^{(\alpha, \beta)}\}_{i = 0}^{\infty}$ for $\mathcal{H}$ by the normalized Jacobi polynomials $P_i^{(\alpha, \beta)}$ (see \cite[Sec.~4.3]{szego1959orthogonal}).
According to the projection theorem (see, e.g. \cite{Luenberger:1997:OVS:524037}) the best fitting (in the least squares sense) approximation of $N$-th order $\hat{x} \in \mathcal{H}$  of $\bar{x}$ always exists unambiguously, and can be calculated by
\begin{equation}
\label{polynomialApprox1}
\hat{x}(\tau) = \sum_{i=0}^{N} \underbrace{\langle \bar{x},P_i^{(\alpha ,\beta)}\rangle}_{c_i} P_i^{(\alpha,\beta)}(\tau), \qquad \tau \in \left[-1,1\right],
\end{equation}
while $c_i$ is considered as the corresponding expansion coefficient. Furthermore, the approximation $\hat{x}$ is orthogonal, i.e. $\langle \hat{x},p \rangle = 0$,  to any polynomial $p \in \bar{\pi}_N \coloneqq \mathrm{span}\{P_i^{(\alpha, \beta)}\}_{i = N +1}^{\infty}$ and is exactly equal to $\bar{x}$, if $\bar{x}\in \pi_N \coloneqq \mathrm{span}\{P_i^{(\alpha, \beta)}\}_{i = 0}^{N} $. 

By applying the transformation \eqref{timeTransform1}, \eqref{timeTransform2} the approximation \eqref{polynomialApprox1} is valid in the time window $\mathcal{I}_{t,T}$ and can be evaluated at any time $t' \in \mathcal{I}_{t,T}$. It is reasonable to choose the evaluation at the time $t$ in order to approximate values at current time. However, by adding a delay $t_d \ge 0$ as a zero $p_{N+1}^{(\alpha,\beta)}$ of the Jacobi polynomial $P_{N+1}^{(\alpha,\beta)}$, the order of the approximation error $\tilde{x} = x-\hat{x}$ can be reduced by one (see \cite{Mboup2009})\label{timedelay}.  For this reason it makes sense to introduce a delay $t_d$ if the delay is justifiable with regard to fault detection.  The delayed polynomial approximation of $x$ based on \eqref{polynomialApprox1} can thus be written as
\begin{align}
\label{polynomialApproxTime}
\hat{x}(t-t_d) &= \sum_{i=0}^{N} \langle x \circ \phi_T,P_i^{(\alpha ,\beta)}\rangle \, (P_i^{(\alpha,\beta)}\circ\phi^{-1}_T(t-t_d)) \nonumber\\
&= \langle x \circ \phi_T,R_{N,t_d}^{(\alpha ,\beta)}\rangle,
\end{align}
with
\begin{equation}
	R_{N,t_d}^{(\alpha ,\beta)}(\tau) = \sum_{i=0}^{N} P_i^{(\alpha ,\beta)}(\tau)\, (P_i^{(\alpha,\beta)}\circ\phi^{-1}_T(t-t_d)).
\end{equation}
The definition of the inner product \eqref{scalar} can be used to represent $\hat{x}(t-t_d)$ by the integral
\begin{equation}
\label{polynomialApproxInt1}
\hat{x}(t-t_d) = \int_{-1}^{1}(x\circ\phi_T(\tau))  g_{N,t_d}(\tau) \mathrm{d}\tau,
\end{equation}
with the kernel 
\begin{equation}
\label{polynomialApproxKernel1}
g_{N,t_d}(\tau) = 	R_{N,t_d}^{(\alpha ,\beta)}(\tau)\, w^{(\alpha ,\beta)}(\tau)\, .
\end{equation}
To evaluate the integral within the original time window $\mathcal{I}_{t,T}$, the substitution $\bar{\tau} = t-\phi_T(\tau)$ is performed. The approximation can therefore be written as
\begin{equation}
\label{polynomialApproxInt2}
\hat{x}(t-t_d) = \int_{0}^{T}x(t-\bar{\tau})  g_{N,t_d}(\bar{\tau}) \mathrm{d}\bar{\tau}  \eqqcolon  \mathcal{P}_{N,t_d}\{x\}(t),
\end{equation}
with the kernel 
\begin{equation}
\label{polynomialApproxKernel2}
g_{N,t_d}(\bar{\tau}) = \dfrac{2}{T} (R_{N,t_d}^{(\alpha ,\beta)}\, w^{(\alpha ,\beta)})\circ \phi^{-1}_T(t-\bar{\tau})\,,
\end{equation}
which is independent of $t$, since $\phi^{-1}_T(t-\bar{\tau}) = 1-\frac{2}{T}\bar{\tau}$. Based on these assumptions, for any $x \in L_2(\left[t-T,t\right]) $ the operator $\mathcal{P}_{N,t_d}\{x\}$ can now be defined as the time-delayed polynomial approximation based on \eqref{polynomialApproxInt2}.

Assume $x \in L_2(\left[t-T,t\right]) \cap C^{k-1}(\left[t-T,t\right])$ and the $k$-th derivative $x^{(k)}$ exists and is Lebesgue integrable.  Since the kernel \eqref{polynomialApproxKernel2} and its $k-1$ derivatives have a compact support due to \eqref{weight} in $\left[0,T\right]$, if $\alpha, \beta \ge k$, the polynomial approximation $x^{(k)}$ can be calculated by successive application of integration by parts
 \begin{align}
 \label{polynomialFIRDer}
 \mathcal{P}_{N,t_d}\{x^{(k)}\}(t) &= \int_{0}^{T} x^{(k)}(t-\tau)  g_{N,t_d}(\tau)\mathrm{d}\tau\nonumber\\
 &=\int_{0}^{T} x(t-\tau)  g_{N,t_d}^{(k)}(\tau)\mathrm{d}\tau \nonumber\\
 &\eqqcolon\mathcal{P}_{N,t_d}^{(k)}\{x\}(t)
 \end{align}
 with the derivative of the kernel given by 
 \begin{align}
 \label{impulseFIRDer}
 g_{N,t_d}^{(k)}(\tau)=(-1)^k\frac{2}{T}(	R_{N,t_d}^{(\alpha ,\beta)}w^{(\alpha,\beta)})^{(k)}\circ \phi^{-1}_T(t-\tau). 
 \end{align} 
 The polynomial approximation of $x^{(k)}$ can thus be calculated by applying of the differentiation approximation operator $\mathcal{P}_{N,t_d}^{(k)}\{\cdot\}$ to $x$.

According to the previously established definition of the polynomial approximation operator, specific properties can be deduced, which are listed below:
\begin{itemize}
	\item \textbf{linearity:} For any two functions $x_1, x_2 \in L_2(\left[t-T,t\right])$, and any constants $c_1, c_2 \in \mathbb{R}$, the equivalence
	  \begin{align}
	  \label{linearity}
	&\mathcal{P}_{N,t_d}\{c_1 x_1 + c_2 x_2 \}(t) \nonumber\\
	&\quad= c_1 \mathcal{P}_{N,t_d}\{x_1 \}(t) + c_2 \mathcal{P}_{N,t_d}\{x_2 \}(t)\, ,
	  \end{align}
	  is valid.
	  \item \textbf{composition:} For a function $x\in L_2(\left[t-T,t\right])$ and a second Lipschitz continuous function $\psi: \mathbb{R}\mapsto \mathbb{R}$ the commutation of the composition of $\psi$ and the polynomial Approximation $\mathcal{P}_{N,t_d}\{x \}(t)$ are approximately equal:
	  \begin{align}
	  \label{composition}
	  &\mathcal{P}_{N,t_d}\{\psi(x) \}(t) \approx \psi\left(\mathcal{P}_{N,t_d}\{x \}(t)\right) = \psi(\hat{x})\, ,
	  \end{align}
	  if $x$ is at least $C^{0}(\left[t-T,t\right])$ and the Lipschitz constant $L$ of $\psi$ is sufficiently small.
	  	  	\item \textbf{reproduction:} For any function $x \in L_2(\left[t-T,t\right])$, with $x\circ \phi_T \in \pi_N$, the polynomial approximation is equal to the approximated function:
	  	\begin{align}
	  		\label{reporduction}
	  		\mathcal{P}_{N,t_d}\{x \}(t) &= x(t-t_d) \nonumber\\
	  		&= \sum_{i=0}^{N}c_i P_i^{(\alpha,\beta)}\circ \phi^{-1}_T(t-t_d).
	  	\end{align}
		The coefficients $c_i$ can be determined by the integral
	  	\begin{align}
	  	\label{reporduction2}
	  	c_i = \int_{0}^{T} x(t-\tau) g_{c_i}(\tau) \mathrm{d}\tau \eqqcolon \mathcal{P}_{N,0,c_i}\{x \}(t),
	  	\end{align}
	  	with the time independent kernel
	  	\begin{align}
	  	g_{c_i}(\tau) = \dfrac{2}{T} (P_{i}^{(\alpha ,\beta)}\, w^{(\alpha ,\beta)})\circ \phi^{-1}_T(t-\tau).
	  	\end{align}
	  Thus, the representation of the function $x$ as a linear combination of Jacobi polynomials to the $N$-th order is equivalent to $x$ if $x \in \pi_N$. The coefficients $c_i$ can be calculated directly from $x$ by \eqref{reporduction2}.
	  	\item \textbf{differentiation:} For each $k$-times differentiable function $x \in L_2(\left[t-T,t\right])$, assuming that the corresponding time derivative of $k$-th degree $x^{(k)}\in L_2(\left[t-T,t\right])$ exists and is Lebesgue integrable, the polynomial approximation of $x^{(k)}$ is equivalent to the differentiation approximation of $x$:
		\begin{align}
			\label{differentiation}
			\mathcal{P}_{N,t_d}\{x^{(k)} \}(t) = \mathcal{P}_{N,t_d}^{(k)}\{x\}(t).
		\end{align}
	 	\item \textbf{partial approximation:} For any two functions $x_1, x_2 \in L_2(\left[t-T,t\right])$ and $x_1 \in \pi_{N^*}$ with $N^* \in \mathbb{N}$, the polynomial approximation of the product $x_1 x_2$ can be defined by
\begin{subequations}
		\label{partialApproximation}
		\begin{align}
		&\mathcal{P}_{N,t_d}\{x_1\, x_2 \}(t) =\mathcal{P}_{N,t_d}\{\mathcal{P}_{N^*,0}\{x_1\}\, x_2\}(t)  \nonumber\\ 
		&=\mathcal{P}_{N,t_d}\{\sum_{i=0}^{N^*}c_i (P_i^{(\alpha,\beta)}\circ\phi^{-1}_T) \, x_2\}(t) \nonumber\\
		&=\sum_{i=0}^{N^*}c_i \mathcal{P}_{N,t_d}\{ (P_i^{(\alpha,\beta)}\circ \phi^{-1}_T)\, x_2\}(t).
		\end{align}
		The components of the polynomial approximation of the product of $x \in L_2(\left[t-T,t\right])$ and the Jacobi polynomial $ P_i^{(\alpha,\beta)}\circ \phi^{-1}_T\in L_2(\left[t-T,t\right])$ can be determined by
		\begin{align}
		\label{modifiedOperator}
		&\mathcal{P}_{N,t_d}\{ (P_i^{(\alpha,\beta)}\circ \phi^{-1}_T)  x\}(t) \nonumber\\
		&= 		\int_{0}^{T}x(t-\tau)  \tilde{g}_{N,t_d,i}(\tau) \mathrm{d}\tau \eqqcolon \tilde{\mathcal{P}}_{N,t_d,i}\{x\}(t),		
		\end{align}
	\end{subequations}
with the kernel 
\begin{equation}
\label{polynomialApproxKernel5}
\tilde{g}_{N,t_d,i}(\tau) = \dfrac{2}{T} (P_i^{(\alpha,\beta)}\,R_{N,t_d}^{(\alpha ,\beta)}\, w^{(\alpha ,\beta)})\circ \phi^{-1}_T(t-\tau)\, .
\end{equation}
Thus, a polynomial approximation of the product $x_1 x_2$ can be obtained by the sum of the products of the coefficients $c_i$ of the approximation $\mathcal{P}_{N^*,0}\{x_1\}$ and the modified  polynomial approximation $\tilde{\mathcal{P}}_{N,t_d,i}\{x_2\}$. It should also be noted that the order $N^*$ of the partial approximation of $x_1$ does not have to be related to the order $N$ of the approximation of the product. It only needs to be large enough to map $x_1$ to the $\mathcal{I}_{t,T}$ interval.  

 Furthermore it is also possible to transfer the differentiation of $x_1$ or $x_2$ to the kernel $\tilde{g}_{N,t_d,i}$ or $g_{c_i}$, respectively, analogous to the differentiation \eqref{polynomialFIRDer} because the property of the compact support has not changed by the multiplication with the Jacobi polynomial of  $i$-th degree.
\end{itemize} 
Based on the presented properties, the polynomial approximation of a complex function with unavailable time derivatives can also be obtained and successively calculated solely by the available signals and the application of the presented integral transformations.
\subsection{Real-time Implementation}
 The polynomial approximation can be calculated from the history of the signal $x$ in the interval $\mathcal{I}_{t,T}$. 
Since the continuous time signal $x$ has to be represented for the realtime-implementation as a sequence of samples $x[k] = x(kT_s)$ of the length $L = \frac{T}{T_s}\in \mathbb{N}$ with the sampling time $T_s \in \mathbb{R}$, a discrete approximation of the integral must be realized. Therefore, the midpoint rule is used to convert the integral to the weighted sum \label{discreteApprox}
\begin{align}
\label{discreteApproximation}
\mathcal{P}_{N,t_d}\{x\left\}[k\right] &=\sum_{j=0}^{L-1}\underbrace{\int_{j\,T_s}^{(j+1)T_s}g_{N,t_d}(\tau) \mathrm{d}\tau}_{w\left[j\right]} \, x\left[k- j\right] \nonumber\\
&= \sum_{j=0}^{L-1} w\left[j\right]  \,x\left[k- j\right],
\end{align}
with the weights $w\left[j\right] \in \mathbb{R},\,j \in \{1,...,L-1\}$. Since the weights $w\left[j\right]$ are constant and therefore independent of the time and the measured signals, they can be determined in advance and do not have to be calculated during evaluation. An analogous procedure to \eqref{discreteApproximation} can also be defined for the differentiation approximation operator \eqref{polynomialFIRDer}, to determine the coefficients \eqref{reporduction2} and the modified approximation operator \eqref{modifiedOperator}, since in all cases it is an integral transformation with a time independent kernel. For this reason, the polynomial approximation can be executed directly on a controller in real-time, by evaluating the individual weighted sums at runtime.

\subsection{Fault detection} 
For fault detection, a residual signal $r$ must be calculated, which fulfils the conditions I and II.  Therefore, equation \eqref{generalMotionA} is rearranged accordingly, to isolate the known signals $q$ and $u$ from the fault $f$. For notational convenience, the dependencies of $q$ and $\dot{q}$ in the matrices $\mathcal{M}(q), \mathcal{C}\left(q,\dot{q}\right),\mathcal{G}\left(q \right) , \mathcal{F}(q,\dot{q}) $ and $ \mathcal{D}(q,\dot{q})$ are not displayed. The generated residual $r_d$, which still depends on $d$, is thus introduced explicitly by 
\begin{align}
r_d &= \mathcal{M}\,\ddot{q} +\mathcal{C}\dot{q}+\mathcal{G} - u =  \Phi(q,\dot{q},\ddot{q},u) = \mathcal{F}f + \mathcal{D}d
\end{align}
with $n_r = n$ and meets the condition I of the RGP-FD  if $d$ is neglected. The second condition is satisfied if $\mathcal{F}(q,\dot{q})$ is non-zero for all $q$ and $\dot{q}$. If disturbances affect the system, it is necessary to decouple the residual from the unknown disturbance $d$. In order to eliminate this dependency, determine the left annihilator $\mathcal{D}^{\perp} \in \mathbb{R}^{n\times n}$ of $\mathcal{D}$, i.e., $\mathcal{D}^\perp\mathcal{D} = 0$, by 
\begin{equation}
\mathcal{D}^{\perp} = I - \mathcal{D}\mathcal{D}^{\dagger}, 
\end{equation}
with the Moore-Penrose generalized inverse $\mathcal{D}^{\dagger} \in \mathbb{R}^{n_d \times n}$  given by 
\begin{align}
\mathcal{D}^{\dagger}  = (\mathcal{D}^\top\mathcal{D})^{-1}\mathcal{D}^\top,
\end{align}
 since $\mathcal{D}$ has full column rank (see, e.g., \cite{penrose_1955}). By premultiplication of the residual with the annihilator $\mathcal{D}^{\perp} \in \mathbb{R}^{n\times n}$, the influence of the disturbance on $r_d$ is eliminated yielding the desired residual
\begin{align}
\label{residualGlobal}
r &= \mathcal{D}^{\perp}\left(\mathcal{M}\,\ddot{q} +\mathcal{C}\dot{q}+\mathcal{G} - u \right)=\mathcal{D}^{\perp}\mathcal{F}f
\end{align}
 according to \eqref{residual1}. It is unequal to the zero vector if any fault $f$ is present, as long as $\mathcal{D}^{\perp}\mathcal{F}_i \ne 0,  \forall i \in \{1,...,n_f\}, \forall q,\dot{q} \in \mathbb{R}^n$. Then $r$ solves the RGP-FD and can be used for fault detection. The only remaining problem is that the determined residual also depends on time derivatives of $q$ and therefore cannot be evaluated. Now the estimates of $\dot{q}$ and $\ddot{q}$ could simply be inserted into the equation \eqref{residualGlobal} to evaluate the residual $r$. However, since $\ddot{q}$ does not have to be continuous, it can happen that the estimation error caused by the polynomial approximation of $\ddot{q}$ and the further multiplication with other terms does not consistently reproduce the residual and thus excites it. For this reason, a polynomial approximation of the whole residual $r$ is performed, by using the polynomial approximation operator $\mathcal{P}_{N,t_d}\{\cdot\}$ instead of replacing $\dot{q}$ and $\ddot{q}$ by the derivative estimates $\hat{\dot{q}}$ and $\hat{\ddot{q}}$ within the residual. The polynomial approximation does not affect the properties of the residual w.r.t. the decoupling of the input $u$ and the disturbance $d$ and can be determined sufficiently well according to the choice of the polynomial degree $N$. By applying the operator to \eqref{residualGlobal} and using its linearity, the polynomial approximation of $r$ can be determined by
\begin{align}
\label{residualGlobalPoly1}
\hat{r}=\mathcal{P}_{N,t_d}\{r\} &= \mathcal{P}_{N,t_d}\{\mathcal{D}^{\perp}\left(\mathcal{M}\,\ddot{q} +\mathcal{C}\dot{q}+\mathcal{G} - u \right)\}\nonumber\\
&=\mathcal{P}_{N,t_d}\{\mathcal{D}^{\perp}\mathcal{M}\,\ddot{q}\} +\mathcal{P}_{N,t_d}\{\mathcal{D}^{\perp}\mathcal{C}\dot{q}\}\nonumber\\
&\quad +\mathcal{P}_{N,t_d}\{\mathcal{D}^{\perp}\mathcal{G}\} - \mathcal{P}_{N,t_d}\{\mathcal{D}^{\perp}u\} \nonumber\\
&=\mathcal{P}_{N,t_d}\{\mathcal{D}^{\perp}\mathcal{F}f\}.
\end{align}
 Every single term can therefore be determined individually by applying the partial approximation \eqref{partialApproximation} and the differentiation approximation operator \eqref{differentiation}. Then, the unknown derivatives of $q$ can be successively eliminated from \eqref{residualGlobalPoly1} and the approximation of the resiudal is given by
 \begin{align}
 \label{residualGlobalPoly2}
 \hat{r}=& \sum_{j=0}^{N^*} \mathcal{P}_{N,0,c_j}\{\mathcal{D}^{\perp}(q,\hat{\dot{q}})\mathcal{M}(q) \} \tilde{\mathcal{P}}^{(2)}_{N,t_d,j}\{q\} \nonumber \\ &\quad +\mathcal{D}^{\perp}(\hat{q},\hat{\dot{q}})\mathcal{C}(\hat{q},\hat{\dot{q}})\hat{\dot{q}} +\mathcal{P}_{N,t_d}\{\mathcal{D}^{\perp}(q,\hat{\dot{q}})\mathcal{G}(q)\}\nonumber\\
 &\quad  - \sum_{j=0}^{N^*} \mathcal{P}_{N,0,c_j}\{\mathcal{D}^{\perp}(q,\hat{\dot{q}}) \} \tilde{\mathcal{P}}_{N,t_d,j}\{u\}
 \end{align}
with $\hat{q} = \mathcal{P}_{N,t_d}\{q\}$, $\hat{\dot{q}} = \mathcal{P}_{N,t_d}^{(1)}\{q\}$. Since \eqref{residualGlobalPoly2} meets criteria I and II of RGP-FD, it can be used for fault detection.

\subsection{Fault identification and estimation}
As shown before, it is possible to use the residual to determine the occurrence of a fault. However, in many cases it is still necessary to know which fault has occurred and to estimate it accordingly. This can be achieved by solving \eqref{residualGlobal} for $f$, which is possible for $\mathrm{rank} ( \mathcal{D}^{\perp}\mathcal{F}) = n_f$ and leads to the fault identification 
\begin{align}
\label{faultDetection}
 f &=(\mathcal{D}^{\perp}\mathcal{F})^{\dagger}\mathcal{D}^{\perp} \left(\mathcal{M}\,\ddot{q} +\mathcal{C}\dot{q}+\mathcal{G}- u \right).
\end{align}
Therein, $(\mathcal{D}^{\perp}\mathcal{F})^{\dagger}$ is the Moore-Penrose generalized inverse of $\mathcal{D}^{\perp}\mathcal{F}$. If $\mathrm{rank} ( \mathcal{D}^{\perp}\mathcal{F})= \tilde{n}_f < n_f$ fault identification can still be partially realized. In this case the solution of \eqref{residualGlobal} is given by
\begin{align}
\label{faultIdentificationExtra1}
f &=(\mathcal{D}^{\perp}\mathcal{F})^{\dagger}\mathcal{D}^{\perp} \left(\mathcal{M}\,\ddot{q} +\mathcal{C}\dot{q}+\mathcal{G}- u \right)+ (\mathcal{D}^{\perp}\mathcal{F})^{\perp}_R l,
\end{align}
in which $(\mathcal{D}^{\perp}\mathcal{F})^{\perp}_R \in \mathbb{R}^{n_f \times (n_f-\tilde{n}_f)}$ is the full column rank right annihilator of $\mathcal{D}^{\perp}\mathcal{F} $, i.e., $\mathcal{D}^{\perp}\mathcal{F} (\mathcal{D}^{\perp}\mathcal{F})^{\perp}_R = 0$ and any vector $l \in \mathbb{R}^{n_f- \tilde{n}_f}$. Hence, $f$ cannot be determined uniquely anymore. However, any linear combination $\tilde{f} = T f$ of $f$ with $T\in \mathbb{R}^{\tilde{n}_f\times n_f}$ and
\begin{align}
T(\mathcal{D}^{\perp}\mathcal{F})^{\perp}_R &=0
\end{align}
can still be computed unambiguously, which follows from premultiplication of \eqref{faultIdentificationExtra1} by $T$. Therein, $T$ is parametrizable by 
\begin{align}
T &= N (\mathcal{D}^{\perp}\mathcal{F}) 
\end{align}
with an arbitrary matrix $N \in \mathbb{R}^{\tilde{n}_f \times n}$ satisfying $\mathrm{rank} N (\mathcal{D}^{\perp}\mathcal{F}) = \tilde{n}_f$ to obtain linear independent faults to be identified. The degrees of freedom in the latter matrix can be utilized to select faults of interest. A solution for $N$ that fulfils $\mathrm{rank} N (\mathcal{D}^{\perp}\mathcal{F}) = \tilde{n}_f$ obviously always exists, since $\mathrm{rank} ( \mathcal{D}^{\perp}\mathcal{F}) = \tilde{n}_f$. Subsequently, however, the case that $\mathrm{rank} ( \mathcal{D}^{\perp}\mathcal{F}) = n_f, \forall q,\dot{q} \in \mathbb{R}^n$ is always considered for fault diagnosis. Analogous to fault detection, by the polynomial approximation the fault estimate can be defined and represented by
\begin{align}
\label{faultDetectionPoly1}
\hat{f} &=\mathcal{P}_{N,t_d}\{\mathcal{K}\left(\mathcal{M}\,\ddot{q} +\mathcal{C}\dot{q}+\mathcal{G} - u \right)\}\nonumber\\
&=\mathcal{P}_{N,t_d}\{\mathcal{K}\mathcal{M}\,\ddot{q}\} +\mathcal{P}_{N,t_d}\{\mathcal{K}\mathcal{C}\dot{q}\}\nonumber\\
&\quad +\mathcal{P}_{N,t_d}\{\mathcal{K}\mathcal{G}\} - \mathcal{P}_{N,t_d}\{\mathcal{K}u\} ,
\end{align}
with $\mathcal{K} = (\mathcal{D}^{\perp}\mathcal{F})^{\dagger}\mathcal{D}^{\perp}$ and the fault estimate $\hat{f} \in \mathbb{R}^{n_f}$. The RGP-FDI can therefore be solved by \eqref{faultDetectionPoly1} if $f$ can be locally approximated within the interval $\mathcal{I}_{t,T}$ by  a polynomial of $N$-th degree.

\section{EXAMPLE}
\tikzstyle{load}   = [ultra thick,-latex]
\tikzstyle{stress} = [-latex]
\tikzstyle{dim}    = [latex-latex]
\tikzstyle{axis}   = [-latex,black]
\pgfmathsetmacro\yangle{100}
\tikzstyle{test} = {[y={({0.5*cos(30)*1cm},{0.5*sin(30)*1cm})},x={(1cm,0cm)},z={(0cm,1cm)}]}
\tikzstyle{isometric}=[x={(0.710cm,-0.410cm)},y={(0cm,0.820cm)},z={(-0.710cm,-0.410cm)}]
\tikzstyle{dimetric} =[x={(0.935cm,-0.118cm)},y={(0cm,0.943cm)},z={(-0.354cm,-0.312cm)}]
\tikzstyle{dimetric2}=[x={(0.935cm,-0.118cm)},z={(0cm,0.943cm)},y={(+0.354cm,+0.312cm)}]
\tikzstyle{trimetric}=[x={(0.926cm,-0.207cm)},y={(0cm,0.837cm)},z={(-0.378cm,-0.507cm)}]
		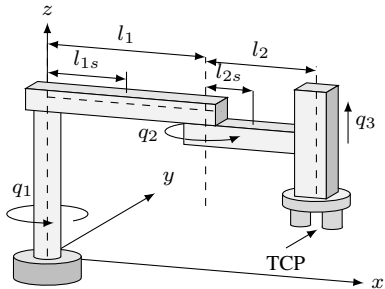
\begin{figure}[h]
			\centering
			\footnotesize
		\begin{tikzpicture}[dimetric2,yscale=0.6,xscale = 0.9]
		\coordinate (O) at (0,0,0);
		\draw[axis] (O) -- ++(5,0,0) node[right] {$x$};
		\draw[axis] (O) -- ++(0,4.5,0) node[above right] {$y$};
		\draw[axis] (O) -- ++(0,0,5.5) node[above] {$z$};
		
		\draw[fill=gray!50] (0,0,-0.5) circle (0.5); 
		\fill[fill=gray!50] (-0.46,-0.2,-0.5) -- (0.46,0.2,-0.5) -- (0.46,0.2,0) -- (-0.46,-0.2,0) -- cycle;
		\draw[fill=gray!20] (O) circle (0.5);
		\draw (0.46,0.2,-0.5) -- ++(0,0,0.5) ;
		\draw (-0.46,-0.2,-0.5) -- ++(0,0,0.5);
		
		 show rotation1-1
		 \draw [shift=(-30:0.6)](0,0,1) arc (-30:270:0.6) node[above left,midway] {$q_1$};

		\draw[fill=gray!10] (O) circle (0.2);
		\fill[fill=gray!10] (-0.175,-0.1,0) -- (0.175,0.1,0) -- ++(0,0,4) -- (-0.175,-0.1,4) -- cycle;
		\draw (-0.175,-0.1,0) -- ++(0,0,4);
		\draw (0.175,0.1,0) -- ++(0,0,4);
		
		 show rotation1-2
		\draw[-latex] [shift=(270:0.6)](0,0,1) arc (270:300:0.6);

		\draw [shift=(20:0.6)](2.5,0,3.255) arc (20:270:0.6) ;

		\draw[fill=gray!20] (2.25,-0.25,3.5) -- (4,-0.25,3.5) -- (4,+0.25,3.5) -- (2.25,+0.25,3.5) -- cycle; 
		\draw[fill=gray!50] (+4.00,-0.25,3) -- (4,+0.25,3) -- (4,+0.25,3.5) -- (+4.00,-0.25,3.5) -- cycle; 
		\draw[fill=gray!10] (2.25,-0.25,3) -- (4,-0.25,3) -- (4,-0.25,3.5) -- (2.25,-0.25,3.5) -- cycle; 

		\draw[fill=gray!20] (-0.25,-0.25,4) -- (2.75,-0.25,4) -- (2.75,+0.25,4) -- (-0.25,+0.25,4) -- cycle; 
		\draw[fill=gray!50] (2.75,-0.25,3.5) -- (2.75,+0.25,3.5) -- (2.75,+0.25,4) -- (2.75,-0.25,4) -- cycle; 
		\draw[fill=gray!10] (-0.25,-0.25,3.5) -- (2.75,-0.25,3.5) -- (2.75,-0.25,4) -- (-0.25,-0.25,4) -- cycle; 

		\draw[fill=gray!30] (4,0,1.25) circle (0.15); 
		\draw[fill=gray!30] (4.5,0,1.25) circle (0.15); 
		\fill[fill=gray!30] (3.862,-0.05,1.25) -- (4.138,0.05,1.25) -- (4.138,0.05,1.75) -- (3.862,-0.05,1.75) -- cycle;
		\fill[fill=gray!30] (4.362,-0.05,1.25) -- (4.638,0.05,1.25) -- (4.638,0.05,1.75) -- (4.362,-0.05,1.75) -- cycle;
		\draw[fill=gray!10] (4,0,1.75) circle (0.15);
		\draw[fill=gray!10] (4.5,0,1.75) circle (0.15);
		\draw (3.862,-0.05,1.25) -- ++(0,0,0.5) ;
		\draw (4.138,0.05,1.25) -- ++(0,0,0.5);
		\draw (4.362,-0.05,1.25) -- ++(0,0,0.5) ;
		\draw (4.638,0.05,1.25) -- ++(0,0,0.5);
		
		\draw[stress] (3.75,0.05,0.65) node[left,below] {TCP} -- ++(0.5,0,0.5) ;
		
		\draw[fill=gray!30] (4.25,0,1.75) circle (0.5); 
		\fill[fill=gray!30] (3.79,-0.2,1.75) -- (4.71,0.2,1.75) -- (4.71,0.2,2) -- (3.79,-0.2,2) -- cycle;
		\draw[fill=gray!10] (4.25,0,2) circle (0.5);
		\draw (4.71,0.2,1.75) -- ++(0,0,0.25) ;
		\draw (3.79,-0.2,1.75) -- ++(0,0,0.25);

		\draw[fill=gray!20] (4,-0.25,4.5) -- (4.5,-0.25,4.5) -- (4.5,+0.25,4.5) -- (4,+0.25,4.5) -- cycle; 
		\draw[fill=gray!50] (4.500,-0.25,2) -- (4.5,+0.25,2) -- (4.5,+0.25,4.5) -- (+4.500,-0.25,4.5) -- cycle; 
		\draw[fill=gray!10] (4,-0.25,2) -- (4.5,-0.25,2) -- (4.5,-0.25,4.5) -- (4,-0.25,4.5) -- cycle; 
		
		\draw[-latex] [shift=(260:0.6)](2.5,0,3.255) arc (260:350:0.6) ;
		\node at (1.6,0,3.04) {$q_2$};

		\draw[dim] (0,0,5.0) -- ++(2.5,0,0) node[midway,above] {$l_1$};
		\draw[dim] (1.25,0,4.325) -- ++(-1.25,0,0) node[midway,above] {$l_{1s}$};
		
		   \draw[-latex] (4.75,0,3.25) -- ++(0,0,1) node[midway,right] {$q_3$};
		\draw[dim] (2.5,0,4.9) -- ++(1.75,0,0) node[midway,above] {$l_2$};
		\draw[dim] (3.25,0,4.3) -- ++(-0.75,0,0) node[midway,above] {$l_{2s}$};

		\draw (3.25,0,3.5) -- ++(0,0,0.8);

		\draw[axis,dashed,-] (0,0,3.75) -- (2.75,0,3.75);
		\draw[axis,dashed,-] (2.5,0,1.5) -- (2.5,0,5);
		\draw[axis,dashed,-] (4.25,0,2) -- (4.25,0,5);

		\coordinate (P) at (1.25,0,4.0);
		\draw (P) -- ++(0,0,0.5);
		\draw[axis,dashed,-] (O) -- (0,0,5);
		\end{tikzpicture} %
		\caption{ Schematic illustration of a SCARA. }
			\label{fig:SCARAPIC}
	\end{figure}
To illustrate the results of this paper, consider a SCARA as an example of a nonlinear mechanical system. This system has three degrees of freedom and consists of three movable segments with the masses $m_1$, $m_2$, $m_3$ and the corresponding inertia values $J_1$, $J_2$ and $J_3$ (see Fig. \ref{fig:SCARAPIC}). The lengths $l_{1s}$ and $l_{2s}$ describe the distance between the center of mass of the first and the second segment from the corresponding rotation axis. For \eqref{generalMotion} the matrices and vectors $\mathcal{M}(q)$, $\mathcal{C}(q,\dot{q})$ and $\mathcal{G}(q)$ are 
\begin{subequations}\label{SCARAMotion}
\begin{align}
	\mathcal{M}(q) &= \begin{bmatrix} \theta_1 + 2 \theta_2 \cos(q_2)&\theta_3 + \theta_2 \cos(q_2)&0\\\theta_3 + \theta_2 \cos(q_2)&\theta_3&0 \\ 0 & 0 & m_{3} \end{bmatrix}\\
	\mathcal{C}(q,\dot{q}) &= \begin{bmatrix} -\dot{q}_2\theta_2 \sin(q_2)& -(\dot{q}_1+ \dot{q}_2 )\theta_2\sin(q_2)&0\\\theta_2 \dot{q}_1\sin(q_2)&0&0 \\ 0 & 0 &0 \end{bmatrix}\\
	\mathcal{G}(q) &= \begin{bmatrix} 0&0 & m_{3}g \end{bmatrix}^\top,
\end{align}
\end{subequations}
with the constant parameters $\theta_1 = {J}_{1}+{J}_{2}+{J}_{3}+l_{1}^2 m_{2}+l_{1}^2 m_{3} +l_{2}^2\,m_{3}+l_{{1s}}^2 m_{1}+l_{{2s}}^2 m_{2},\theta_2 = l_{1}\,l_{2}\,m_{3}\,+\,l_{1}\,l_{{2s}}\,m_{2} $ and $\theta_3 = m_{3}\,l_{2}^2+m_{2}\,l_{2s}^2+{J}_{2}+{J}_{3}$.

In addition, actuator faults $f_1$ and $f_2$ are considered for the first two joints, and a force acting in $y$-$z$-direction at the TCP with an unknown value is regarded as the disturbance $d$. The matrices $\mathcal{F}(q,\dot{q})$ and $\mathcal{D}(q,\dot{q})$ are thus given by
\begin{subequations}
\begin{align}\label{SCARAfault}
\mathcal{F}(q,\dot{q}) &= \begin{bmatrix} 1&0&0\\0&1&0 \end{bmatrix}^\top\\
\mathcal{D}(q,\dot{q}) &= \begin{bmatrix} l_{2}\cos\left(q_{1}+q_{2}\right)+l_{1}\cos\left(q_{1}\right)\\ l_{2}\cos\left(q_{1}+q_{2}\right)\\ 1 \end{bmatrix},
\end{align}
\end{subequations}
 with $\mathrm{rank}\, \mathcal{D} = 1 = n_d$. The entire dynamic of the SCARA can thus be described as a nonlinear equation of motion. Since $\mathrm{rank}\,\mathcal{D}^\perp\mathcal{F}= 2 = n_f$ the fault can be determined according to \eqref{faultDetection} by  
\begin{subequations}\label{faultEstimationSCARA}
\begin{align}
\label{faultEstimationSCARA1}
f_1 &=\theta_1 \ddot{q}_1 +  2 \theta_2 \cos(q_2) \ddot{q}_1 +\theta_3 \ddot{q}_2 - g l_2 m_3 \cos(q_1+q_2)  \nonumber\\
&\quad  + l_1 \cos(q_1) u_3 -  l_2 m_3\cos(q_1+q_2) \ddot{q}_3 +   \theta_2 \cos(q_2) \ddot{q}_2\nonumber\\
&\quad -  l_1 m_3\cos(q_1) \ddot{q}_3- g l_1 m_3\cos(q_1) - \dot{q}^2_2 \theta_2 \sin(q_1)\nonumber\\
&\quad -  2 \dot{q}_1 \dot{q}_2 \theta_2 \sin(q_1)  + l_2\cos(q_1 + q_2) u_3 - u_1\\
f_2 &= \theta_3 \ddot{q}_1+ \theta_2 \cos(q_2)\ddot{q}_1 + \theta_3 \ddot{q}_2  - u_2+\dot{q}_1^2 \theta_2 \sin(q_1)\nonumber\\
&\quad+ l_2\cos(q_1+q_2) u_3 - g l_2 m_3 \cos(q_1+q_2) \nonumber\\
&\quad   - m_3 l_2\cos(q_1+q_2) \ddot{q}_3 \label{faultEstimationSCARA2}  .
\end{align}
\end{subequations}
In order to obtain a fault estimate $\hat{f}$, it is necessary to apply the polynomial approximation to \eqref{faultEstimationSCARA} and then use the operator's introduced properties \eqref{linearity} to \eqref{partialApproximation} in such a way that the dependencies on time derivatives of the input or output signals are eliminated. The estimated faults $\hat{f}_1$ and $\hat{f}_2$ then become
\begin{subequations}	\label{faultEstimationSCARAges2}
	\begin{align}
\label{faultEstimationSCARA3}
\hat{f}_1 &=\theta_1 \hat{\ddot{q}}_1 +  2 \theta_2 \mathcal{P}_{N,t_d}\{\cos(q_2)\ddot{q}_1\}  +\theta_3 \hat{\ddot{q}}_2  - \hat{\dot{q}}^2_2 \theta_2 \sin(\hat{q}_1) \nonumber \\
&\quad - g l_2 m_3 \cos(\hat{q}_1+\hat{q}_2) + l_1\mathcal{P}_{N,t_d}\{\cos(q_1)u_3\}   \nonumber\\
&\quad -  l_2 m_3  \mathcal{P}_{N,t_d}\{\cos(q_1+ q_2)\ddot{q}_3\} +   \theta_2 \mathcal{P}_{N,t_d}\{\cos(q_2)\ddot{q}_2\}\nonumber \\
&\quad -  l_1 m_3\mathcal{P}_{N,t_d}\{\cos(q_1)\ddot{q}_3\}- g l_1 m_3\cos(\hat{q}_1)\nonumber\\
&\quad -  2 \hat{\dot{q}}_1 \hat{\dot{q}}_2 \theta_2 \sin(\hat{q}_1)  + l_2 \mathcal{P}_{N,t_d}\{\cos(q_1+ q_2)u_3\} - \hat{u}_1\\
\hat{f}_2 &= \theta_3 \hat{\ddot{q}}_1+ \theta_2 \mathcal{P}_{N,t_d}\{\cos(q_2)\ddot{q}_1\} + \theta_3 \hat{\ddot{q}}_2  - \hat{u}_2\nonumber\\
&\quad +\hat{\dot{q}}_1^2 \theta_2 \sin(\hat{q}_1) + l_2 \mathcal{P}_{N,t_d}\{\cos(q_1+ q_2)u_3\} \nonumber\\
&\quad- g l_2 m_3 \cos(\hat{q}_1+\hat{q}_2) - m_3 l_2 \mathcal{P}_{N,t_d}\{\cos(q_1+ q_2)\ddot{q}_3\},\label{faultEstimationSCARA4}
\end{align}
\end{subequations}
with the estimated values $\hat{q}_1,\hat{q}_2,\hat{\dot{q}}_1,\hat{\dot{q}}_2,\hat{\ddot{q}}_1,\hat{\ddot{q}}_2, \hat{u}_1, \hat{u}_2$ and the partial approximation of the remaining products. The resulting fault estimate $\hat{f}$ fulfils criteria I and II of RGP-FDI and is solely dependent on the available measurement signals.

For the simulation the physical parameters are defined according to the Table \ref{table_example}. An exemplary pick and place movement with a quadratic spline interpolation between four positions was generated as the desired trajectory. In order to realize the trajectory tracking for the robot, a PI state controller with an additional feedforward control to compensate the nonlinearity was implemented. All eigenvalues of the state controller were set to $\lambda = -10 $. The parameters for the polynomial approximation have been chosen to $\alpha= \beta = 3 $ and $N = 1$, i.e., first order Jacobi polynomials are imployed. In order to improve the approximation accuracy the delay $t_d$ was selected as the zero of the Jacobi polynomial of second order at $t_d = (L \, T_s)/3$ (see Sec. \ref{timedelay}) and correspondingly for discrete-time implementation of the simulation the sampling time of $T_s = 0.005 \si{\second}$ and $L = 20$ were set  (see Sec. \ref{discreteApprox}). Furthermore, a measurement noise $\bar{\omega}_y$ and a process noise $\bar{\omega}_u$ for the variables $q$ and $u$, respectively, were added for the simulation to verify the robustness of the method against noise. The disturbance $d$ can be described by $d = \sigma(t-0.5\si{s})(10 \si{N}+2\si{N} \sin(2t))$, whereby the faults $f_1$ and $f_2$ jump to the value $f_{1,\infty} = f_{2,\infty} = 10 \si{Nm}$ at the time $t_{f1}= 1\si{s}$ and $t_{f2}= 3\si{s}$, respectively.
{\renewcommand{\arraystretch}{1.3}%
\begin{table}[h]
\caption{simulation parameters}
\label{table_example}
\begin{center}
\begin{tabular}{|l|l||l|l|}
\hline
parameter & value&parameter & value\\
\hline
\hline
$m_1$ & $10\si{\kilogram}$&$J_3$ & $0.005\si{\kilogram} \si{\metre}^2$\\
\hline
$m_2$ & $5\si{\kilogram}$&$l_1$ & $0.325\si{\metre}$\\
\hline
$m_3$ & $2.35\si{\kilogram}$&$l_{1s}$ & $0.1625\si{\metre}$\\
\hline
$J_1$ & $0.088\si{\kilogram} \si{\metre}^2$&$l_2$ & $0.275\si{\metre}$\\
\hline
$J_2$ & $0.0315\si{\kilogram} \si{\metre}^2$&$l_{2s}$ & $0.1375\si{\metre}$\\
\hline
$g$ & $9.81 \si{\metre} \si{\second}^{-2}$&&\\
\hline

\end{tabular}
\end{center}
\end{table}}

\begin{figure}
	\centering
	\begin{tikzpicture}
	\begin{axis}[name=plot4,
	footnotesize,
	scaled ticks=false,
	width=8.5cm,height=3cm,
	xmin=0,xmax=5,
	ymin=-0.0006,ymax=0.0006,
	ytick={-0.0005,0,...,0.0005},
	xtick={0,0.5,...,5},
	ytick align=inside,
	xtick align=inside,
	xticklabels=\empty,
	legend pos = south east
    ]
	\addplot [blue] table[x=t, y = noise_y1]{data.dat};\label{pny1}	
\addlegendentry{$\bar{\omega}_{y1}$}
\addplot [purple] table[x=t, y = noise_y2]{data.dat};\label{pny2}	
\addlegendentry{$\bar{\omega}_{y2}$}
\addplot [green] table[x=t, y = noise_y3]{data.dat};\label{pny3}	
\addlegendentry{$\bar{\omega}_{y3}$}
	\end{axis}
	\begin{axis}[
	name = plot5,
	footnotesize,
	at={($(plot4.south)-(0,0.1cm)$)},
	anchor=north,
	width=8.5cm,height=3cm,
	xmin=0,xmax=5,
	ymin=-1.6,ymax=1.6,
	ytick={-1,0,...,1},
	xtick={0,0.5,...,5},
	ytick align=inside,
	xtick align=inside,
	xticklabels=\empty,
	legend pos = south east
	]
	\addplot [blue] table[x=t, y = noise_u1]{data.dat};\label{pnu1}	
	\addlegendentry{$\bar{\omega}_{u1}$}
	\addplot [purple] table[x=t, y = noise_u2]{data.dat};\label{pnu2}	
\addlegendentry{$\bar{\omega}_{u2}$}
	\addplot [green] table[x=t, y = noise_u3]{data.dat};\label{pnu3}	
\addlegendentry{$\bar{\omega}_{u3}$}
	\end{axis}
		\begin{axis}[name=plot6,
	footnotesize,
	width=8.5cm,height=3cm,
	at={($(plot4.south)-(0,1.6cm)$)},
	anchor=north,
	xmin=0,xmax=5,
	ymin=-0.5,ymax=18,
	ytick={0,5,...,19},
	xtick={0,0.5,...,5},
	ytick align=inside,
	xtick align=inside,
	xticklabels=\empty,
	legend pos = south east
	]
	\addplot [black,thick] table[x=t, y = d]{data.dat};\label{p3}	
	\addlegendentry{$d$}
	\end{axis}
	\begin{axis}[
	name = plot7,
	footnotesize,
	at={($(plot4.south)-(0,3.1cm)$)},
	anchor=north,
	width=8.5cm,height=3cm,
	xmin=0,xmax=5,
	ymin=-3,ymax=14,
	ytick={0,5,...,13},
	xtick={0,0.5,...,5},
	ytick align=inside,
	xtick align=inside,
	xticklabels=\empty,
	legend pos = south east
	]
	\addplot [black,dashed,thick] table[x=t, y = f1]{data.dat};\label{p4}	
	\addlegendentry{$f_{1}$}
	\addplot [red,thick] table[x=t, y = fest1]{data.dat};	\label{p5}
	\addlegendentry{$\hat{f}_{1}$}
	
	\end{axis}
	\begin{axis}[name = plot8,
	footnotesize,
	width=8.5cm,height=3cm,
	xmin=0,xmax=5,
	ymin=-3,ymax=13,
	ytick={0,5,...,14},
	at={($(plot4.south)-(0,4.6cm)$)},
	anchor = north,
	legend pos = south east
	]
	\addplot [black,dashed,thick] table[x=t, y = f2]{data.dat};\label{p6}	
\addlegendentry{$f_{2}$}
\addplot [red,thick] table[x=t, y = fest2]{data.dat};	\label{p7}
\addlegendentry{$\hat{f}_{2}$}
	\end{axis}
	\end{tikzpicture}	
	\caption{Applied disturbance $d$, the obtained faults $f_1$, $f_2$ and their estimates $\hat{f}_1$, $\hat{f}_2$ in the presence of measurement and process noise.}
		\label{faults}
\end{figure}
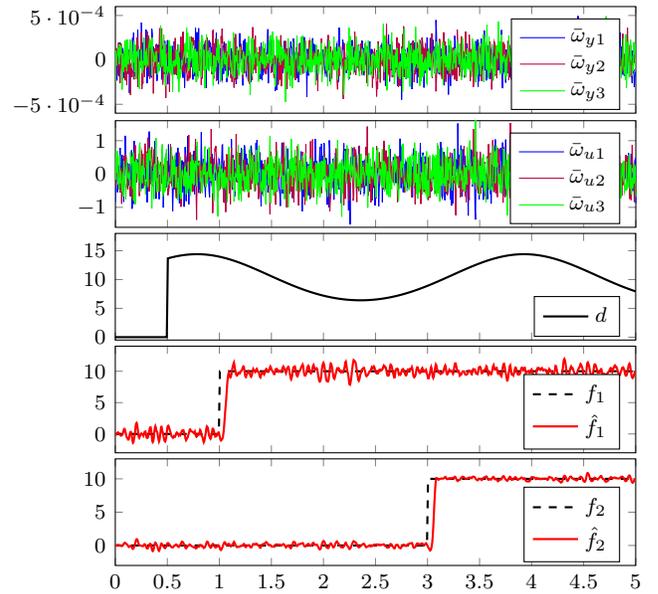

The simulation results in Fig. \ref{movements} and Fig. \ref{faults} show the trajectories of $q$, $u$, $f$, and $d$ for the exemplary pick and place movement. It becomes obvious that the reconstructed values $\hat{f}_1$ and $\hat{f}_2$ of the faults $f_1$ and $f_2$ are each independent of the input variable $u$ and the disturbance $d$.  The reconstruction of the fault thus depends solely on the parameters of the polynomial approximation. The simulation results show that the faults were both reconstructed with a delay of $0.033 \si{\second}$, which corresponds to the set delay $t_d$ of the polynomial approximation. For the reconstruction it does not matter whether the respective measured variables or the disturbance can be approximated polynomially within the interval  $\mathcal{I}_{t,T}$ or not. The only requirement is that the fault $f$ can be sufficiently well approximated within the interval  $\mathcal{I}_{t,T}$ by a polynomial of $N$-th degree.

\begin{figure}[h]
	\centering
	
	\begin{tikzpicture}
\begin{axis}[name=plot1,
footnotesize,
width=9cm,
height=3cm,
xmin=0,xmax=5,
ymin=0.15,ymax=0.4,
ytick={0,0.1,...,0.5},xtick={0,0.5,...,5},
ytick align=inside,
xtick align=inside,
 xticklabels=\empty,
legend pos = south east]
\addplot [black,thick] table[x=t, y = y1]{data.dat};\label{p1}	
\addlegendentry{$q_1$}
\addplot [blue,dashed,thick] table[x=t, y = y1s]{data.dat};	\label{p2}
\addlegendentry{$q_{1s}$}
\end{axis}
\begin{axis}[name = plot2,
footnotesize,
at={($(plot1.south)-(0,0.1cm)$)},
anchor=north,
width=9cm,
height=3cm,
xmin=0,xmax=5,
ymin=-0.1,ymax=0.8,
ytick={0,0.2,...,0.7},
xtick={0,0.5,...,5},
 xticklabels=\empty,
ytick align=inside,xtick align=inside,
legend pos = south east
]
\addplot [black,thick] table[x=t, y = y2]{data.dat};
\addlegendentry{$q_2$}	
\addplot [blue,dashed,thick] table[x=t, y = y2s]{data.dat};
\addlegendentry{$q_{2s}$}	
\end{axis}
\begin{axis}[name = plot3,
footnotesize,
width=9cm,
height=3cm,
xmin=0,xmax=5, 
ymin=0.1,ymax=0.45,
ytick={0.2,0.3,...,0.5},
at={($(plot1.south)-(0,1.6cm)$)},
anchor = north,
xticklabels=\empty,
legend pos = north east
]
\addplot [black,thick] table[x=t, y = y3]{data.dat};
\addlegendentry{$q_{3}$}	
\addplot [blue,dashed,thick] table[x=t, y = y3s]{data.dat};
\addlegendentry{$q_{3s}$}	
\end{axis}
\begin{axis}[name = plot9,footnotesize,width=9cm,height=3.0cm,xmin=0,xmax=5, ymin=-25,ymax=30,ytick={-20,-10,...,20},at={($(plot1.south)-(0,3.1cm)$)},anchor = north,
legend pos = north east
]
\addplot [blue,thick] table[x=t, y = u1]{data.dat};\label{pu1}
\addlegendentry{$u_{1}$}
\addplot [purple,thick] table[x=t, y = u2]{data.dat};\label{pu2}
\addlegendentry{$u_{2}$}
\addplot [green,thick] table[x=t, y = u3]{data.dat};\label{pu3}
\addlegendentry{$u_{3}$}

\end{axis}
\end{tikzpicture}	
\caption{Simulation results are the position values and the references as well as the corresponding control inputs for a pick and place movement of the SCARA.}
\label{movements}
\end{figure}
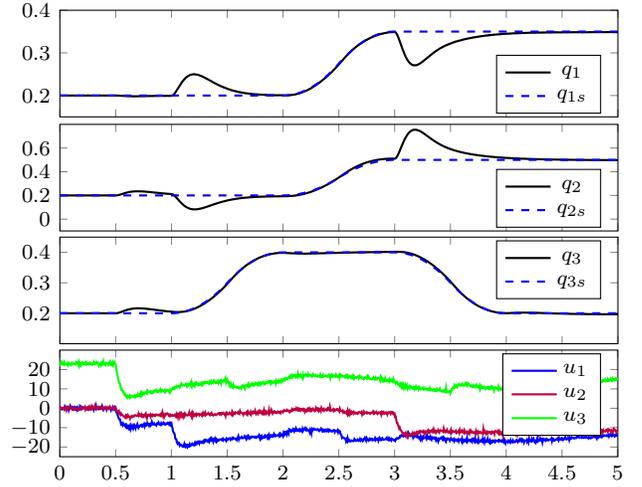

\section{CONCLUDING REMARKS}
As discussed in this paper, faults in nonlinear mechanical systems can be detected or identified independently by the introduced polynomial approximation independent of the disturbance $d$ and the form of the input or output signals. It should be noted that the presented method can be extended to other influences such as friction and elastic couplings in the joints. In addition to actuator faults, other faults, such as sensor and parameter faults, can also be detected and identified analogous to the presented method. In further research work, a more general class of nonlinear systems will also be investigated.

\addtolength{\textheight}{-12cm}   



%

\section*{ACKNOWLEDGMENT}
The authors kindly express their gratitude to the industrial
research partner Siemens AG, Digital Industry Division
Erlangen for funding and supporting this project.


\bibliographystyle{IEEEtran}
\bibliography{IEEEabrv,literatur}

\end{document}